\documentclass[leqno]{amsart}
\usepackage{amsmath}
\usepackage{amssymb}
\usepackage{amsthm}
\usepackage{enumerate}
\usepackage[mathscr]{eucal}
\theoremstyle{plain}
\usepackage{tikz}
\usepackage{faktor}
\newtheorem{theorem}{Theorem}[section]
\newtheorem{lemma}[theorem]{Lemma}

\theoremstyle{definition}
\newtheorem{definition}[theorem]{Definition}
\newtheorem{remark}[theorem]{Remark}

\newtheorem{example}[theorem]{Example}

\newtheorem{cor}[theorem]{Corollary}
\theoremstyle{remark}
\usepackage{hyperref}
\hypersetup{
    linkcolor=blue,
    filecolor=red,      
    urlcolor=magenta,
}

\begin{document}

\newcommand{\Z}{\ensuremath{\mathcal{Z}}}
\newcommand{\tr}{\ensuremath{{}^t\!}}
\newcommand{\tra}{\ensuremath{{}^t}}
\newcommand{\J}{\mathsf{J}}
\newcommand{\M}{\mathsf{M}}
\newcommand{\W}{\mathsf{W}}
\newcommand{\D}{\mathsf{D}}
\newcommand{\V}{\mathsf{V}}
\newcommand{\X}{\mathbb{X}}
\newcommand{\C}{\mathbb{C}}
\newcommand{\R}{\mathfrak{Re}}
\newcommand{\G}{\mathsf{G}}
\newcommand{\U}{\mathsf{U}}
\newcommand{\NU}{\mathsf{NU}}
\newcommand{\Inn}{\textup{Inn}}
\newcommand{\Out}{\textup{Out}}

\newcommand{\diag}{\textup{diag}}
\newcommand{\gal}{\textup{Gal}}
\newcommand{\cross}{\times}
\newcommand{\tensor}{\otimes}
\newcommand{\secref}[1]{Section~\ref{#1}}
\newcommand{\thmref}[1]{Theorem~\ref{#1}}
\newcommand{\lemref}[1]{Lemma~\ref{#1}}
\newcommand{\remref}[1]{Remark~\ref{#1}}
\newcommand{\propref}[1]{Proposition~\ref{#1}}
\newcommand{\corref}[1]{Corollary~\ref{#1}}
\newcommand{\eqnref}[1]{~{\textrm(\ref{#1})}}
\newcommand{\tabref}[1]{Table~(\ref{#1})}
\newcommand{\figref}[1]{Figure~(\ref{#1})}
\newcommand{\defref}[1]{Definbition~(\ref{#1})}
\makeatletter
\def\imod#1{\allowbreak\mkern10mu({\operator@font mod}\,\,#1)}
\makeatother

\theoremstyle{plain}
\theoremstyle{definition}

\numberwithin{equation}{section}
\newcommand{\comment}[1]{\marginpar{#1}}
\newcommand{\ignore}[1]{}
\newcommand{\term}[1]{\textit{#1}}
\newcommand{\mynote}[1]{}
\def \GL {\mathrm{GL}}
\def \SL {\mathrm{SL}}
\def \Sp {\mathrm{Sp}}

\title[smoothness in operator spaces induced by numerical radius]{Numerical radius and a notion of smoothness in the space of bounded linear operators}

                \newcommand{\acr}{\newline\indent}

\author{Saikat Roy, Debmalya Sain}

\address{(Roy)~Department of Mathematics, National Institute of Durgapur, West Bengal, India}
\email{saikatroy.cu@gmail.com}

\address{(Sain)~Department of Mathematics, Indian Institute of Science, Bengaluru 560012, Karnataka, India}
\email{saindebmalya@gmail.com}


\thanks{The research of Dr. Debmalya Sain is sponsored by SERB N-PDF Post-doctoral fellowship, under the mentorship of Professor Apoorva Khare. Dr. Sain feels elated to acknowledge the friendship of Dr. Milind Jha, his childhood friend and an accomplished doctor. The research of Mr. Saikat Roy is supported by CSIR MHRD in form of Senior Research Fellowship under the supervision of Professor Satya Bagchi. }

\begin{abstract}
We observe that the classical notion of numerical radius gives rise to a notion of smoothness in the space of bounded linear operators on certain Banach spaces, whenever the numerical radius is a norm. We demonstrate an important class of \emph{real} Banach space $ \mathbb{X} $ for which the numerical radius defines a norm on $ \mathbb{L}(\mathbb{X}), $ the space of all bounded linear operators on $ \mathbb{X} $. We characterize Birkhoff-James orthogonality in the space of bounded linear operators on a finite-dimensional Banach space, endowed with the numerical radius norm. Some examples are also discussed to illustrate the geometric differences between the numerical radius norm and the usual operator norm, from the viewpoint of operator smoothness.
\end{abstract}
\maketitle

\section{Introduction}

The main objective of the present article is to study smoothness in the space of bounded linear operators on a Banach space, induced by the numerical radius. The study of smoothness in the space of bounded linear operators on a Banach space, \emph{with respect to the usual operator norm}, is a classical area of research in geometry of Banach spaces \cite{A,H,HE,HO,KY,PSG,R1,R2,R3,SPMR,W}. The space of bounded linear operators on a Banach space, endowed with the numerical radius norm, need not be isometrically isomorphic to the space of bounded linear operators on the same Banach space, endowed with the usual operator norm, in general. Therefore, it is expected to have differences in geometric structures in the space of bounded linear operators on a Banach space, equipped with these two different norms. The current work explores the said differences from the point of view of smoothness.\\

The symbol $\mathbb{X}$ signifies a Banach space over the field $\mathbb{F}$, where $\mathbb{F}=\mathbb{R}$ or $\C$. Unless otherwise mentioned, we work with both real and complex Banach spaces. Given any $\lambda\in \C$, let $\R~\lambda$ denote the real part of $\lambda$. For any subset $\D$ of $\mathbb{F}$, let $\mathsf{CO}(\D)$ denote the convex hull of $\D.$ It is immediate that if $\D$ is a compact subset of $\mathbb{F}$ then $\mathsf{CO}(\D)$ is also a compact subset of $\mathbb{F}$.\\

Let $B_\X$ and $S_\X$ denote the closed unit ball and the unit sphere of $\X$, respectively. Let $ext(B_\X)$ denote the collection of all extreme points of the closed unit ball $B_\X$. We denote the zero vector of any vector space by $\theta$, other than the scalar field $\mathbb{F}$. Let $\mathbb{L}(\X)$ $(\mathbb{K}(\X))$ denote the collection of all bounded (compact) linear operators on $\X$ endowed with the usual operator norm. We use the symbol $\M_T$ to denote the norm attainment set of a bounded linear operator $T\in \mathbb{L}(\X)$, i.e., $\M_T=\{x\in S_\X:\|Tx\|=\|T\|\}.$ Given any $x,y\in \X$, we say that $x$ is Birkhoff-James orthogonal \cite{B} to $y$, written as $x\perp_B y$, if
$\|x+\lambda y\|\geq \|x\|$ for all scalars $\lambda \in \mathbb{F}.$ Let $\X^*$ denote the topological dual of $\X$. The collection of all support functionals at a non-zero $x \in \X $ is denoted by $J(x)$ and is defined by:
\begin{align*}
J(x) : = \left\{ x^*\in S_{\X^*} : x^*(x) = \|x\|\right\}.
\end{align*}
By the Hahn-Banach Theorem, the collection $J(x)$ is non-empty. It is well known \cite{J,Ja} that  $x\perp_B y$ if and only if there exists $x^*\in J(x)$ such that $x^*(y)=0$. The element $x$ is said to be a smooth point in $\X$, if $J(x)$ is singleton. Equivalently, $x$ is a smooth point in $\X$ if and only if for any $y_1,y_2\in \X$ with $x\perp_B y_1$ and $x\perp_B y_2$ imply that $x\perp_B (y_1+y_2)$, i.e., Birkhoff-James orthogonality is right additive at $x$. The space $\X$ is called smooth if every non-zero element of $\X$ is smooth.\\

Given any $T \in \mathbb{L}(\X), $ the numerical range of $T$ is defined by:
\begin{align*}
\W(T) :=\left\{ x^*(Tx) :  (x,x^*)\in \J \right\},~\mathrm{where}~\mathsf{J} : =  \left\{(x,x^*)\in \mathbb{X}\times {\mathbb{X}^*}: x\in S_\mathbb{X}, x^*\in J(x)\right\}.
\end{align*}
The numerical radius of the linear operator $T$ is defined by:
$$\|T\|_w =\sup\{|\lambda|:\lambda \in \W(T)\}.$$
It is well known that whenever $\mathbb{F}=\C$, the numerical radius defines a norm on $\mathbb{L}(\X)$. However, $\|\cdot\|_w$ need not be a norm on $\mathbb{L}(\X)$, if $\mathbb{F}=\mathbb{R}$. Throughout the text, we will only consider those Banach spaces for which $\|\cdot\|_w$ defines a norm in $\mathbb{L}(\X)$. The space of bounded linear operators on $\X$ endowed with the numerical radius norm is denoted by $(\mathbb{L}(\X))_w$. For a detailed study on numerical range of operators and their possible applications, we refer the readers to \cite{GK, KMMPQ, MM, MMR}.\\

Birkhoff-James orthogonality is an important tool in the study of smoothness of elements in a given Banach space. Indeed, using the Birkhoff-James orthogonality of operators \cite{S}, a complete characterization of smoothness in $\mathbb{K}(\X)$ was provided in \cite{PSG}, for a real reflexive Banach space $\X$. The geometry of $\mathbb{L}(\X)$ is heavily dependent on the norm attainment sets of its members \cite{BS, Sa, Sai, SP}. In particular, the norm attainment set of a linear operator $T$ plays a pivotal role in determining the smoothness of $T$ in $\mathbb{L}(\X)$. A characterization of smoothness without any restriction on $\M_T$ was provided in \cite{SPMR}. The above studies motivate  us to explore the concept of {\emph{numerical radius smoothness}} in $(\mathbb{L}(\X))_w$.

\begin{definition}\label{nu-orthogonality}
Let $\X$ be a Banach space and let $T,A\in (\mathbb{L}(\mathbb{X}))_w$. We say that $T$ is numerical radius Birkhoff-James orthogonal to $A$, written as $T\perp_B^w A$ if
$$\|T+\lambda A\|_w\geq \|T\|_w \qquad \forall~ \lambda\in \mathbb{F}.$$
\end{definition}

\begin{definition}\label{nu-smooth}
Let $\X$ be a Banach space and let $T\in (\mathbb{L}(\mathbb{X}))_w$ be non-zero. We say that $T$ is {\emph{nu-smooth}} (the abbreviated form of numerical radius smooth), if
$$T\perp_B^w A,~T\perp_B^w B~\mathrm{imply~that}~T\perp_B^w (A+B)\qquad \forall~ A,B\in (\mathbb{L}(\X))_w.$$
\end{definition}

It is easy to see that given any $T\in (\mathbb{L}(\X))_w$,
\begin{align}\label{Polygonal reformulation}
\|T\|_w =\sup\left\{|x^*(Tx)|: (x,x^*)\in \left(ext(B_\X)\times ext(B_{\X^*})\right)\bigcap \J\right\}.    
\end{align}
The above formulation is particularly advantageous whenever $\mathbb{X}$ is a finite-dimensional real polyhedral Banach space, i.e., $ext(B_\X)$ is finite. Note that a finite-dimensional real Banach space $\X$ is polyhedral if and only if $\X^*$ is polyhedral and a member $x^*$ of $\X^*$ is an extreme point of $B_{\X^*}$ if and only if $x^*$ is the unique supporting functional corresponding to a facet of $B_\X$ \cite[Lemma 2.1]{SPBB}. For any non-zero element $T\in (\mathbb{L}(\X))_w,$ let $\M_{\W(T)}$ denote the the numerical radius attainment set of $T$, i.e.,
\begin{align*}
\M_{\W(T)} := \{(x,x^*)\in \J : x^*(Tx) = \sigma,~|\sigma|= \|T\|_w \}. 
\end{align*}
Note that $\M_{\W(T)}$ is non-empty, whenever $\X$ is finite-dimensional. The collection of support functional at $T$ is defined by:
\begin{align*}
J_\W(T): = \{ f\ : \left(\mathbb{L}(\X)\right)_w \to \mathbb{F} : f~\mathrm{is~linear},\|f\|=1, f(T) = \|T\|_w \}.
\end{align*}
It follows from the James characterization (\cite[Theorem 2.1]{Ja}) that $T\perp^w_B A$, for some $A\in \left(\mathbb{L}(\X)\right)_w$ if and only if there exists $f\in J_\W(T)$ such that $A\in \ker f.$ Also, for any non-zero $T\in \left(\mathbb{L}(\X)\right)_w$, $T$ is nu-smooth if and only if $J_\W(T)$ is singleton.\\

After this introductory part, the present article is demarcated in two sections. In Section \ref{Section 3}, we provide a large class of real finite-dimensional Banach spaces, where numerical radius defines a norm on the concerned spaces of linear operators. We acquire a characterization of smoothness in $(\mathbb{L}(\X))_w$, both in the finite-dimensional and the infinite-dimensional cases, in Section \ref{Section 4}. In due course of our development, we also obtain a necessary and sufficient condition for Birkhoff-James orthogonality in $(\mathbb{L}(\X))_w$ whenever $\X$ is finite-dimensional. Some examples have been discussed to show that smoothness in $\mathbb{L}(\X)$ and smoothness in $(\mathbb{L}(\X))_w$ are not equivalent.

\section{Numerical radius and polyhedral Banach spaces}\label{Section 3}

In \cite{MM}, the authors studied numerical radius of some two-dimensional real polygonal Banach spaces and obtained explicit formulae for the numerical index of the said spaces. In particular, this illustrates that numerical radius is a norm on $\mathbb{L}(\X)$ for specific two-dimensional real polygonal Banach space $\X$. The following result shows that numerical radius is a norm on $\mathbb{L}(\X)$ for any finite-dimensional polyhedral Banach space $\X.$

\begin{theorem}\label{Polygonal norm}
Let $\mathbb{X}$ be an $n$-dimensional real polyhedral Banach space. Then $\|\cdot\|_w$ is a norm on $\mathbb{L}(\X)$. \end{theorem}

\begin{proof}
Note that $\mathbb{X}$ is isometrically isomorphic to $ (\mathbb{R}^n, \|.\|), $ for some norm $ \|.\| $ on $ \mathbb{R}^n. $ Let $T\in \mathbb{L}(\mathbb{X})$ with $\|T\|_w=0.$ It is enough to show that $T$ is the zero operator on $\mathbb{X}$. We now complete the proof of the theorem in the following two steps :\\

\noindent Step I:  Let $F_1,F_2, \dots, F_k$ denote the facets of $B_\X$ and let $f_i$ denote the unique support functional corresponding to the facet $F_i$ for each $1\leq i \leq k$. Obviously,
$$S_\X = \left\{x\in \X : \max\left\{f_i(x)\right\}_{i=1}^k=1\right\}.$$

In this step, we show that for any extreme point $x_0$ of $B_\mathbb{X}$, there exist at least $n$ number of distinct facets of $B_\mathbb{X}$ that contain $x_0$. Let $x_0=(z_1,z_2, \dots , z_n).$ Without loss of generality, let $F_1, F_2, \dots , F_r$ be the only facets of $B_\X$ that contain $x_0$. Suppose on the contrary that $r<n.$ We now consider the following sets:
\begin{align*}
& P_1 := \{ f_i : f_i(x_0)=1,~1\leq i \leq k \},\\
& P_2 := \{ f_i : f_i(x_0)=-1,~1\leq i \leq k \},\\
& P_3 := \{ f_i : f_i(x_0)\neq \pm 1,~1\leq i \leq k \}.
\end{align*}
Evidently, $P_1=\{f_i\}_{i=1}^r$, $P_2=\{-f_i: f_i\in P_1\}$ and $\max\left\{ |f_i(x_0)| : f_i\in P_3\right\} < 1- \delta,$ for some $\delta \in (0,1)$. Since $r<n$, we can find some non-zero $u_0\in \X$ such that the following two conditions are satisfied: 
\begin{align*}
(1)~ u_0\in \bigcap\limits_{i=1}^r\ker f_i,~
(2)~ \max\left\{ |f_i(u_0)| : f_i\in P_3\right\} < \dfrac{\delta}{2}.
\end{align*}
Now, we consider the vectors $x_0+u_0$ and $x_0-u_0$ in $\X.$ Then we have that
\begin{align*}
& f_i(x_0\pm u_0) =1 \qquad \forall f_i\in P_1, \\
& f_i(x_0\pm u_0) =-1 \qquad \forall f_i\in P_2,\\
\max & \left\{ |f_i(x_0\pm u_0)| : f_i\in P_3 \right\} <1-\dfrac{\delta}{2}.
\end{align*}
As a result, $x_0\pm u_0\in S_\X$. This proves that $x_0$ is not an extreme point of $B_\mathbb{X},$ which is a contradiction.\\

\noindent Step II: In this step we show that $Tx_0=\theta$ for any $x_0\in ext(B_\mathbb{X}).$ Applying Step - I, we can find $n$ number of distinct facets of $B_\mathbb{X}$ that contain $x_0$, say $F_1,F_2, \dots, F_n$. Let for each $i\in \{1,2, \dots, n\}$, $(u^i_n)$ be a sequence of smooth points in $F_i$ converging to $x_0$. Since $\|T\|_w=0,$ $f_i(Tu^i_n)=0$ for each $i\in \{1,2, \dots , n \}.$ Therefore, $Tu^i_n\in \ker f_i$, for each $i\in \{1,2, \dots , n \}.$ Due to the continuity of $T$, $Tu^i_n\to Tx_0$ and $Tx_0\in \ker f_i$ for each $i\in \{1,2, \dots , n \}.$ Consequently, 
$$Tx_0\in \bigcap\limits_{i=1}^n \ker f_i = \{\theta\}.$$
Thus, $Tx_0=\theta.$ Since $x_0\in ext(B_\mathbb{X})$ was chosen arbitrarily, it follows that $T$ maps every extreme point of $B_\mathbb{X}$ to $\theta,$ as expected.\\

However, since $\M_T\cap ext(B_\X)\neq \emptyset$, $T$ must be the zero operator on $\X$. This completes the proof.
\end{proof}

\begin{remark}\label{remark about polyhedral}
Let $\mathbb{X}$ be an $n$-dimensional real polyhedral Banach space. It is worth mentioning that the above result can also be proved using \cite[Theorem 2.2]{SPBB}. In the said theorem, the authors showed the numerical index of $\X$ is non-zero, assuming that given any extreme point $x_0$ of $B_\X$, there is $n$ number of distinct facets of $B_\X$ that meet at $x_0$. This essentially shows that numerical radius is a norm on $\mathbb{L}(\X)$, whenever any extreme point $x_0$ of $B_\X$ is contained in $n$ number of distinct facets of $B_\X$. However, step I of the above result shows that the said assumption is redundant. Therefore, step I, in combination with \cite[Theorem 2.2]{SPBB} ultimately proves Theorem \ref{Polygonal norm}.
\end{remark}

It is evident that the notion of support functional plays a vital role in the previous theorem. This is also true as far as our next result is concerned. It was proved in \cite{EK} that numerical radius defines a norm on $\mathbb{L}(\ell_p)$, where $1\leq p < \infty$. In the next theorem, we prove the same, using support functionals. We believe that the following proof is simpler in comparison to the proof given in \cite{EK}.

\begin{theorem}\label{lp norm}
Numerical radius is a norm on $\mathbb{L}(\ell_p)$, where $1\leq p < \infty$; $p\neq 2$.
\end{theorem}

\begin{proof}
It is enough to show that $\|T\|_w=0$ implies that $T$ is the zero operator on $\ell_p.$ A Schauder basis for $\ell_p$ is $(e_k)$, where $e_k=(\delta_{kj})$ has $1$ in the $k$-th coordinate and zeros otherwise. For each $j\in \mathbb{N},$ let
$$T(e_j)=\sum\limits_{k\in \mathbb{N}}a_{kj}e_k, \qquad a_{kj}\in \mathbb{R}.$$
For each $x=(x_1,x_2, \dots)\in S_{\ell_p}$, let $x^*:\ell_p\to \mathbb{R}$ be defined by:
$$x^*(y)=\begin{cases}\sum\limits_{k\in \mathbb{N}}y_k~sgn (x_k)|x_k|^{p-1}; &\mathrm{if}~1<p<\infty,\\
\sum\limits_{k\in \mathbb{N}}y_k~sgn (x_k); &\mathrm{if}~p=1,
\end{cases}$$
for all $y=(y_1, y_2, \dots)\in \ell_p$. Obviously, $x^*\in J(x)$ for $p=1$, and therefore, $(x,x^*)\in \J$. Whenever $1<p<\infty$, using H\"{o}lders inequality, it is also not difficult to see that $(x,x^*)\in \J$. It follows from the hypothesis of the theorem that $x^*(Tx)=0$ for all $(x,x^*)\in \J$. Therefore, we obtain
$$e_j^*(Te_j)=a_{jj}=0\qquad \forall~ j\in \mathbb{N}.$$
On the other hand, for any $r,s\in \mathbb{N}$ with $r\neq s,$ we have that
$$\left(\frac{\alpha e_r+\beta e_s}{\|\alpha e_r+\beta e_s\|}\right)^*\left(T\left(\frac{\alpha e_r+\beta e_s}{\|\alpha e_r+\beta e_s\|}\right)\right)=0 \qquad \forall~ \alpha, \beta\in \mathbb{R}\setminus \{0\}.$$
On simplification, we get
$$\alpha\beta(a_{rs}|\alpha|^{p-2}+a_{sr}|\beta|^{p-2})=0\qquad \forall~ \alpha, \beta\in \mathbb{R}\setminus \{0\}.$$
This shows that $a_{rs}=a_{sr}=0.$ Since $r,s\in \mathbb{N}$ was chosen arbitrarily, it follows that $T(e_j)=0$ for all $j\in \mathbb{N}.$ Therefore, $T$ is the zero operator on $\ell^p$, as desired.
\end{proof}

Since $\ell_\infty^n$ is polyhedral, combining Theorem \ref{Polygonal norm} and Theorem \ref{lp norm}, we have the following corollary. We omit the proof, as it is immediate.

\begin{cor}\label{General lp space}
Let $\X=\ell_p^n$, where $n$ is a natural number and $1\leq p\leq \infty;$ $p\neq 2$. Then $\|\cdot\|_w$ is a norm on $\mathbb{L}(\X)$.
\end{cor}

\section{Smoothness induced by the numerical radius}\label{Section 4}

We devote this section to study nu-smoothness of bounded linear operators, which is the integral theme of the present article. We start by characterizing nu-smoothness of a bounded linear operator on any Banach space $\X$.

\begin{theorem}\label{Characterization of smoothness}
Let $\X$ be a Banach space and let $T\in {\left(\mathbb{L}(\X)\right)_w}$ be non-zero. Then the following conditions are equivalent:\\
(i) $T$ is nu-smooth.\\
(ii) $T\perp^w_B A$ for $A\in \left(\mathbb{L}(\X)\right)_w$ implies that for any sequence $\left((x_n,x_n^*)\right)\subseteq \J$ with the property that 
\begin{align}\label{infinite dimension}
lim~x_n^*(Tx_n) \to \sigma,\quad |\sigma|=\|T\|_w,
\end{align}
every sub-sequential limit of the sequence $\left(x_n^*(Ax_n)\right)$ is zero.
\end{theorem}

\begin{proof}
$(i)\implies (ii):$ Suppose on the contrary that there exists a sequence $\left((x_n,x_n^*)\right)\subseteq \J$ satisfying (\ref{infinite dimension}) such that
$$\lim x_{n_k}^*\left(Ax_{n_k}\right)=r \neq 0,$$
for some sub-sequence $\left(x_{n_k}^*\left(Ax_{n_k}\right)\right)$ of $\left(x_n^*(Ax_n)\right)$. Let $B\in \left(\mathbb{L}(\X)\right)_w$ be defined by $B=T-\dfrac{\sigma}{r}A.$ Then we obtain
$$\lim x_{n_k}^*\left(Bx_{n_k}\right)=\lim x_{n_k}^*\left(\left(T-\dfrac{\sigma}{r}A\right)x_{n_k}\right)=\lim x_{n_k}^*\left(Tx_{n_k}\right)- \dfrac{\sigma}{r}\lim x_{n_k}^*\left(Ax_{n_k}\right)= 0.$$
This leads us to conclude that
\begin{align*}
\|T+\lambda B \|_w \geq \lim \left|x_{n_k}^*(T+\lambda B)(x_{n_k})\right|= \lim  \left|x_{n_k}^*(Tx_{n_k})+\lambda x_{n_k}^*(Bx_{n_k})\right| = \|T\|_w,
\end{align*}
for all scalars $\lambda$. In other words, $T\perp^w_B B$. Since $\perp^w_B$ is homogeneous and $T$ is nu-smooth, we get that $T\perp^w_B \left(\dfrac{\sigma}{r}A+B\right) =T$, which is a contradiction.\\

\noindent $(ii) \implies (i):$ Suppose that $T\perp^w_B A$ and $T\perp^w_B B$ for some non-zero $A,B \in \left(\mathbb{L}(\X)\right)_w.$ Consider any sequence $\left((x_n,x_n^*)\right)\subseteq \J$ that satisfies the condition (\ref{infinite dimension}). It now follows from the hypothesis of the theorem that we can find monotonically increasing sequence of natural numbers, say $(n_k)$, such that $$\lim x_{n_k}^*(Ax_{n_k}) =\lim x_{n_k}^*(Bx_{n_k})=0.$$
Therefore, $\lim x_{n_k}^*\left((A+ B)(x_{n_k})\right)=0.$ Now, for every scalar $\lambda,$ we have that
$$\|T+\lambda(A+B)\|_w\geq \lim \left|x_{n_k}^*\left(T+\lambda(A+ B)\right)(x_{n_k})\right|=\|T\|_w.$$
In other words, $T\perp^w_B (A+B).$ Thus, $T$ is nu-smooth and the proof follows.
\end{proof}

An interesting query on this context is whether the above characterization takes any special form if $\X$ is finite-dimensional. An extra advantage in assuming $\X$ to be finite-dimensional is that we now have $\M_{\W(T)}\neq \emptyset$ for any $T\in (\mathbb{L}(\X))_w$. Therefore, we can expect $\M_{\W(T)}$ to play an important role in determining the nu-smoothness of $T$. To explore the said connection, we first prove a lemma that is particularly helpful in our further developments.

\begin{lemma}\label{convex}
Let $\X$ be a finite-dimensional Banach space and let $T,A\in \left(\mathbb{L}(\X)\right)_w$ be non-zero. The set $\D$ defined by 
\begin{align}\label{definition of D}
\D := \left\{ \overline{x^*(Tx)}x^*(Ax) : (x,x^*)\in \M_{\W(T)}\right\},
\end{align}
is a compact subset of $\mathbb{F}.$
\end{lemma}

\begin{proof}
It is trivial to see that $\D$ is bounded. Therefore, to show that $ \D $ is compact, it is sufficient to show that $\D$ is closed. Assume that $(\mu_n)$ is a sequence in $\D$ with $\mu_n \to \mu_0.$ Obviously, for each $n$ 
$$\mu_n = \overline{x_n^*(Tx_n)}x_n^*(Ax_n),~\mathrm{where}~(x_n,x^*_n)\in \M_{\W(T)}.$$ 
Passing through a suitable sub-sequence if necessary, we may assume that $x_n\to x_0$ and $x_n^*\to x_0^*$, as $n\to \infty$, where $x_0\in S_\X$ and $x_0^*\in S_{\X^*}$. Observe that
\begin{align*}
|x^*_n(Tx_n) - x_0^*(Tx_0)| & = |x_n^*(Tx_n) - x_n^*(Tx_0) + x_n^*(Tx_0) - x_0^*(Tx_0)|\\
& \leq |x_n^*(Tx_n - Tx_0)| + |(x_n^* - x_0^*)(Tx_0)|\\
& \leq \|Tx_n-Tx_0\| + \|x_n^* - x_0^* \|\|Tx_0\|.
\end{align*}
Since $T$ is continuous, $Tx_n \to Tx_0,$ as $n\to \infty.$ Thus, $x_n^*(Tx_n)\to x_0^*(Tx_0)$, as $n\to \infty.$ Evidently, $x_0^*(Tx_0) = \sigma$ for some $\sigma\in \mathbb{F}$ with $|\sigma|=\|T\|_w.$ Similar argument shows that $x_n^*(Ax_n)\to x_0^*(Ax_0)$ and $x^*_n(x_n)\to x_0^*(x_0)=1.$ This proves that
$$(x_0,x_0^*)\in \M_{\W(T)}~\mathrm{and}~\lim \overline{x_n^*(Tx_n)}x_n^*(Ax_n) = \mu_0 =\overline{x_0^*(Tx_0)}x_0^*(Ax_0).$$
Consequently, $\mu_0 \in \D.$ Thus, $\D$ is closed and this completes the proof of the lemma.
\end{proof}

The following theorem completely characterizes Birkhoff-James orthogonality in $(\mathbb{L}(\X))_w$ for a finite-dimensional Banach space $\X$.

\begin{theorem}\label{Characterization of orthogonality}
Let $\X$ be a finite-dimensional Banach space and let $T,A\in \left(\mathbb{L}(\X)\right)_w$ be non-zero. Then the following conditions are equivalent:\\
$(i)$ $T\perp_B^w A.$\\
$ (ii) $ $0\in \mathsf{CO}(\D),$ where $\D$ is the subset of $\mathbb{F}$ defined by (\ref{definition of D}).
\end{theorem}

\begin{proof}
$ (i) \implies (ii): $ Since numerical radius Birkhoff-James orthogonality is homogeneous, without loss of generality, we may assume that $\|A\|_w=1$. Suppose on the contrary that $0\notin \mathsf{CO}(\D).$ Since $\mathsf{CO}(\D)$ is a compact convex subset of $\mathbb{F}$ (Lemma \ref{convex}), rotating $\mathsf{CO}(\D)$ suitably if necessary, we may and do assume that $\R~ d > 0$ for all $d\in \mathsf{CO}(\D).$ Moreover, due to the compactness of $\D$, we can find $r\in (0, \frac{1}{2})$ such that $\R~ d > r$ for all $d\in \D.$ In other words,
\begin{align}\label{Compactness of D}
\R~ \overline{x^*(Tx)}x^*(Ax) > r\qquad \forall~ (x,x^*)\in \M_{\W(T)}.
\end{align}
Next, we define
\begin{align*}
\G : = \left\{ (x,x^*) \in \J : \R~ \overline{x^*(Tx)}x^*(Ax) \leq \dfrac{r}{2} \right\}.
\end{align*}
We claim that 
\begin{align*}
\sup\{|x^*(Tx)|:(x,x^*)\in \G\} < \|T\|_w-2\varepsilon\qquad \mathrm{for~some~} \varepsilon \in (0,\frac{1}{2}).
\end{align*}
It follows from (\ref{Compactness of D}) and the definition of $\G$ that $\G \cap \M_{\W(T)} = \emptyset.$ Suppose that $\left((x_n,x_n^*)\right)\subseteq \G $ with $\lim |x_n^*(Tx_n)| = \|T\|_w.$ Without loss of generality, we may assume that $x_n\to x_0$ and $x_n^*\to x_0^*$, as $n\to \infty$, where $x_0\in S_\X$ and $x_0^*\in S_{\X^*}$. Applying the similar techniques as in the proof of Lemma \ref{convex}, it can be shown that 
$$(1)~x_n^*(Tx_n)\to x_0^*(Tx_0)\quad (2)~x_n^*(Ax_n) \to x_0^*(Ax_0)\quad (3)~x_0^*(x_0)=1.$$
Since $|x_0^*(Tx_0)|=\|T\|_w$, we have that $(x_0,x_0^*)\in \M_{\W(T)}$. Therefore,  
$$\R~ \overline{x_n^*(Tx_n)}x_n^*(Ax_n) \to \R~ \overline{x_0^*(Tx_0)}x_0^*(Ax_0) > r\qquad (\mathrm{using}(\ref{Compactness of D})).$$ 
However, this is a contradiction, as $\left((x_n,x_n^*)\right)\subseteq \G$.
Thus, $\sup\{|x^*(Tx)|:(x,x^*)\in \G\} < \|T\|_w-2\varepsilon~\mathrm{for~some~} \varepsilon \in (0,\frac{1}{2}).$

Choose $0 < \lambda < \min \left\{ \varepsilon, r\right\}. $ Now, for any $(x,x^*)\in \G$
\begin{align*}
\left|x^*\left(Tx- {\lambda} Ax\right)\right|  & \leq  \left|x^*(Tx)\right| + \left|{\lambda} x^*(Ax)\right|\\
& < \|T\|_w - 2 \varepsilon + \lambda \\
& < \|T\|_w - \varepsilon.
\end{align*}
Also, for any $(x,x^*)\in \J \setminus \G$
\begin{align*}
\left|x^*\left(Tx- {\lambda} Ax\right)\right|^2 & = x^*\left(Tx- {\lambda} Ax\right) \overline{x^*\left(Tx- {\lambda} Ax\right)}\\
& \leq \|T\|_w^2 + \lambda^2 - 2\lambda \R~ \overline{x^*(Tx)}x^*(Ax)\\
& \leq \|T\|_w^2 + \lambda^2 - \lambda r.
\end{align*}
Since $\lambda^2 - \lambda r < 0$, we get
\begin{align*}
\left\| T- {\lambda} A\right\|_w = \sup \left\{ \left|x^*\left(Tx- {\lambda} Ax\right)\right| : (x,x^*)\in \J \right\} < \|T\|_w. 
\end{align*}
This is a contradiction to the fact that $T\perp^w_B A.$ Therefore, $0\in \mathsf{CO}(\D)$, as desired.\\

\noindent $ (ii) \implies (i): $ Since $0\in \mathsf{CO}(\D),$ applying Carath\'{e}odory Theorem we can find $t_j\in [0,1]$ and $(x_j,x_j^*) \in \M_{\W(T)},$ $j=1,2,3;$ such that
\begin{align}\label{Constructing linear functional}
\sum \limits_{j=1}^3 t_j=1~\mathrm{and}~ \sum \limits_{j=1}^3 t_j \overline{x_j^*(Tx_j)}x_j^*(Ax_j) = 0. 
\end{align}
Let $\rho : \left(\mathbb{L}(\X)\right)_w \to \mathbb{F}$ be defined by 
\begin{align*}
\rho(B) = \dfrac{1}{\|T\|_w}\sum \limits_{j=1}^3 t_j \overline{x_j^*(Tx_j)}x_j^*(Bx_j)\qquad \forall~ B \in \left(\mathbb{L}(\X)\right)_w.
\end{align*}
Clearly, $\rho(A)=0.$ Also, note that for any $B \in \left(\mathbb{L}(\X)\right)_w$, 
\begin{align*}
|\rho(B)| = \left|\dfrac{1}{\|T\|_w}\sum \limits_{j=1}^3 t_j \overline{x_j^*(Tx_j)}x_j^*(Bx_j)\right| \leq \dfrac{1}{\|T\|_w}\sum \limits_{j=1}^3 t_j| \overline{x_j^*(Tx_j)}||x_j^*(Bx_j)| \leq \|B\|_w,  
\end{align*}  
and 
\begin{align*}
\rho(T)=  \dfrac{1}{\|T\|_w}\sum \limits_{j=1}^3 t_j \overline{x_j^*(Tx_j)}x_j^*(Tx_j) =  \dfrac{1}{\|T\|_w}\sum \limits_{j=1}^3 t_j \|T\|_w^2 = \|T\|_w. 
\end{align*}
This shows that $\rho \in J_\W(T).$ Therefore, $T \perp^w_B A$ and the proof follows.
\end{proof}

Whenever $\M_{\W(T)} = \left\{ \left(\mu x_0, \overline{\mu}x_0^*\right) : |\mu|=1, (x_0,x_0^*)\in \J\right\}$, for some fixed $x_0\in S_{\mathbb{X}}$ and $x_0^*\in S_{\mathbb{X}^*}$, we have the following corollary:

\begin{cor}\label{One point norm attainment}
Let $\X$ be a finite-dimensional Banach space and let $T,A\in \left(\mathbb{L}(\X)\right)_w$ be non-zero with $\M_{\W(T)} = \left\{ \left(\mu x_0, \overline{\mu}x_0^*\right) : |\mu|=1, (x_0,x_0^*)\in \J\right\}.$ Then $T \perp^w_B A$ if and only if $x_0^*(Ax_0)=0.$ 
\end{cor}
\begin{proof}
It follows from Theorem \ref{Characterization of orthogonality} that 
\begin{align*}
T \perp^w_B A ~\Leftrightarrow ~ 0 \in \mathsf{CO}\left(\left\{ \overline{\overline{\mu}x_0^*(T\mu x_0)}\overline{\mu}x_0^*(A\mu x_0) : |\mu|=1, (x_0,x_0^*)\in \J\right\}\right). 
\end{align*}
Clearly, 
\begin{align*}
\overline{\overline{\mu}x_0^*(T\mu x_0)}\overline{\mu}x_0^*(A\mu x_0)= \overline{x_0^*(T x_0)}x_0^*(A x_0).   
\end{align*}
As a result, $T \perp^w_B A$ if and only if $x_0^*(Ax_0)=0.$ This completes the proof.
\end{proof}

Finally, we characterize nu-smoothness in $(\mathbb{L}(\X))_w$, for a finite-dimensional Banach space $\X.$

\begin{theorem}\label{Characterization of NU-smooth}
Let $\X$ be a finite-dimensional complex Banach space and let $T\in \left(\mathbb{L}(\X)\right)_w$ be non-zero. Then the following conditions are equivalent : \\
$(i)$ $T$ is nu-smooth.\\
$(ii)$ $\M_{\W(T)} = \left\{ \left(\mu x_0, \overline{\mu}x_0^*\right) : |\mu|=1, (x_0,x_0^*)\in \J\right\}.$
\end{theorem}

\begin{proof}
$(i)\implies (ii):$ Suppose on the contrary that there exists $(y_0,y_0^*)\in \M_{\W(T)}$ such that $(y_0,y_0^*)\neq (\mu x_0,\overline{\mu}x_0^*)$ for any unimodular constant $\mu.$ We now complete the proof of the theorem by considering the following two cases: \\

\noindent Case I: Let $y_0 = \sigma_0 x_0$ for some unimodular constant $\sigma_0.$ We claim that $\ker x_0^* \neq \ker y_0^*.$ Indeed, if  $\ker x_0^* = \ker y_0^*, $ then $y_0^* = \alpha_0 x_0^*$ for some unimodular scalar $\alpha_0.$ Since $y_0^*(\sigma_0 x_0) = 1$, we get $\alpha_0 \sigma_0 = 1$, which is true if and only if $\alpha_0 = \overline{\sigma_0}.$ However, this proves that $y_0^*=\overline{\sigma_0}x_0^*,$ which is a contradiction, since $(y_0,y_0^*)\neq (\mu x_0,\overline{\mu}x_0^*)$ for any unimodular constant $\mu.$ Therefore, $\ker x_0^* \neq \ker y_0^*$, as we have claimed.\\

Next, we consider $x_1,x_2\in \X$ such that $x_1\in \ker x_0^*\setminus \ker y_0^*$ and $x_2\in \ker y_0^*\setminus \ker x_0^*.$ Observe that for any $z\in \X$, there exist a unique scalar $\alpha_z$ and a unique vector $h_z\in \ker x_0^*$ such that
$$z = \alpha_z x_0+ h_z.$$
Now, we define $A_1,A_2: \X\to \X$ by 
\begin{align*}
A_1(z) = \alpha_z x_1~\mathrm{and}~A_2(z)= \alpha_z x_2 \qquad \forall~ z\in \X.
\end{align*}
Clearly, $A_1,A_2 \in \left(\mathbb{L}(\X)\right)_w.$ Note that $A_1(x_0)=x_1$ and $A_2(y_0)=\sigma_0 A_2(x_0)=\sigma_0 x_2.$ Since $(x_0,x_0^*)$, $(y_0,y_0^*)\in \M_{\W(T)}$, we have that
$$x_0^*(Tx_0) = \sigma_1\|T\|_w~\mathrm{and}~y_0^*(Ty_0) = \sigma_2\|T\|_w,$$ 
for some unimodular constant $\sigma_1$, $\sigma_2.$
We define $\rho, \tau : \left(\mathbb{L}(\X)\right)_w\to \mathbb{F}$ by
\begin{align*}
\rho(B) = x_0^*(Bx_0)~\mathrm{and}~\tau(B) = y_0^*(By_0)\qquad \forall~ B\in \left(\mathbb{L}(\X)\right)_w.   
\end{align*}
Observe that $\rho$ and $\tau$ are linear and the linear functional $\overline{\sigma_1}\rho : \left(\mathbb{L}(\X)\right)_w\to \mathbb{F}$ satisfies the following:
\begin{align*}
&~ (i)~ |\overline{\sigma_1}\rho(B)| = |\rho(B)| \leq \|B\|_w\qquad \forall~ B\in \left(\mathbb{L}(\X)\right)_w,\\
&~(ii)~ \overline{\sigma_1}\rho(T) = \|T\|_w,\\
&~(iii) \overline{\sigma_1}\rho(A_1)=\overline{\sigma_1}x_0^*(A_1x_0)=\overline{\sigma_1}x_0^*(x_1)=0,\\
&~(iv)~ \overline{\sigma_1}\rho(A_2)=\overline{\sigma_1}x_0^*(A_2x_0)=\overline{\sigma_1}x_0^*(x_2)\neq 0.\\
\end{align*}
Therefore, we get $\overline{\sigma_1}\rho\in J_\W(T)$, $A_1\in \ker \overline{\sigma_1}\rho$ and $A_2\notin \ker \overline{\sigma_1}\rho.$ Similar arguments show that $\overline{\sigma_2}\tau\in J_\W(T)$, $A_2\in \ker \overline{\sigma_2}\tau$ and $A_1\notin \ker \overline{\sigma_2}\tau.$ Thus, $\overline{\sigma_1}\rho$, $\overline{\sigma_2}\tau$ are distinct members of $J_\W(T)$. As a result, $T$ is not nu-smooth, which is a contradiction.\\

\noindent Case II: Let $y_0\neq \sigma x_0$ for any unimodular scalar $\sigma.$ Let $z_0\in \X$ be such that
\begin{align}\label{z_0}
z_0 = x_0^*(y_0)x_0-y_0.   
\end{align}
Evidently, $z_0$ is non-zero, as otherwise, $x_0^*(y_0)x_0=y_0$ and $|x_0^*(y_0)|=1.$ Observe that $z_0\in \ker x_0^*.$ Consider any $z_0^*\in J(z_0).$ Evidently, for any $z\in \X$ there exist unique scalars $\alpha_z, \beta_z$ and $h_z\in \ker x_0^*\cap \ker z_0^*$ that 
\begin{align}\label{unique expression I}
z = \alpha_z x_0 + \beta_z z_0 + h_z.
\end{align}
Plugging the expression of $z_0$ (see (\ref{z_0})) into (\ref{unique expression I}), we get
\begin{align*}
z = (\alpha_z + \beta_z x_0^*(y_0)) x_0 + (-\beta_z) y_0 + h_z.
\end{align*}
As a result, for every $z\in \X$, there exist $\gamma_z, \zeta_z\in \mathbb{F}$ and $h_z\in \ker x_0^*\cap \ker z_0^*$ such that 
\begin{align*}
z = \gamma_z x_0 + \zeta_z y_0 + h_z.
\end{align*}
Now, we define $T_1,T_2: \X\to \X$ by 
\begin{align*}
T_1(z) = \gamma_z x_0~\mathrm{and}~T_2(z)= \zeta_z y_0\qquad \forall~ z\in \X.
\end{align*}
Clearly, $T_1,T_2 \in \left(\mathbb{L}(\X)\right)_w.$ Note that $T_1(y_0)=T_2(x_0)=\theta.$ Moreover, since $(x_0,x_0^*),(y_0,y_0^*)$ are contained in $\M_{\W(T)}$, we have that
$$x_0^*(Tx_0) = \sigma_1\|T\|_w~\mathrm{and}~y_0^*(Ty_0) = \sigma_2\|T\|_w,$$
for some unimodular scalars   $\sigma_1$, $\sigma_2$. We define $\psi, \eta : \left(\mathbb{L}(\X)\right)_w\to \mathbb{F}$ by
\begin{align*}
\psi(B) = x_0^*(Bx_0)~\mathrm{and}~\eta(B) = y_0^*(By_0)\qquad \forall~ B\in \left(\mathbb{L}(\X)\right)_w.   
\end{align*} 
Consider the linear functionals $\overline{\sigma_1}\psi, \overline{\sigma_2} \eta : \left(\mathbb{L}(\X)\right)_w\to \mathbb{F}$. Then using analogous techniques as in Case I, we can show the following:\\
\begin{align*}
& (i)~ \overline{\sigma_1}\psi, \overline{\sigma_2}\eta \in J_\W(T),\\
& (ii)~ T_2\in \ker \overline{\sigma_1}\psi,~ T_1\notin \ker \overline{\sigma_1}\psi,\\
& (iii)~ T_1\in \ker \overline{\sigma_2}\eta,~ T_2\notin \ker \overline{\sigma_2}\eta.\\
\end{align*}
Thus, $\overline{\sigma_1}\psi$, $\overline{\sigma_2}\eta$ are distinct members of $J_\W(T)$. As a result, $T$ is not nu-smooth, which is a contradiction.\\

\noindent $(ii) \implies (i):$ Suppose that $T \perp^w_BU_1$ and $T \perp^w_BU_2,$ for some non-zero $U_1,U_2 \in \left(\mathbb{L}(\X)\right)_w.$ Then it follows from Corollary \ref{One point norm attainment} that $x_0^*(U_1x_0)=0$ and $x_0^*(U_2x_0)=0.$ Therefore, we have that $x_0^*((U_1+U_2)x_0)=0.$ As a result, $T \perp^w_B(U_1 +U_2)$. This proves that $T$ is nu-smooth and thereby establishes the theorem completely.
\end{proof}

Equipped with the above characterization, we are now in a position to explore the geometrical dissimilarities between $\mathbb{L}(\X)$ and $(\mathbb{L}(\X))_w$ from the perspective of smoothness. We deduce through examples that smoothness in $(\mathbb{L}(\X))_w$ and $\mathbb{L}(\X)$ are not equivalent. Our examples involve finite-dimensional real polyhedral Banach spaces. Therefore, we explicitly mention the above theorem in the real case separately for the convenience of the readers. Also, in this context, it is worth mentioning that $\|\cdot\|_w$ defines a norm in the space of bounded linear operators on a finite-dimensional real polyhedral Banach space [Theorem \ref{Polygonal norm}, Section \ref{Section 3}].

\begin{theorem}\label{Characterization of NU-smooth in R}
Let $\X$ be a finite-dimensional real Banach space and let $T\in \left(\mathbb{L}(\X)\right)_w$ be non-zero. Then the following conditions are equivalent:\\
$(i)$ $T$ is nu-smooth.\\
$(ii)$ $\M_{\W(T)} = \left\{ \left(a x_0, a x_0^*\right) : a \in \{-1,1\}, (x_0,x_0^*)\in \J\right\}.$
\end{theorem}

To serve our purpose, we also state the following result, which characterizes smoothness in the space of compact linear operators on a real reflexive Banach space, endowed with the usual operator norm. Given a Banach space $\X$ and a normed linear space $\mathbb{Y}$, let $\mathbb{K}(\X, \mathbb{Y})$ denote the space of all compact linear operators from $\X$ to $\mathbb{Y}.$

\begin{theorem}\cite[Theorem 4.1 and Theorem 4.2]{PSG}\label{PSG-4.1,4.2}
Let $\X$ be a real reflexive Banach space and let $\mathbb{Y}$ be a real normed space. Then $T \in \mathbb{K}(\X, \mathbb{Y})$ is smooth if and only if $T$ attains norm at a unique (upto scalar multiplication) vector $x_0$ (say) of $S_\X$ and $Tx_0$ is a smooth point.
\end{theorem}

\begin{example}\label{NU smooth but not operator smooth}
Let $\mathbb{Z}=\X\oplus_\infty \mathbb{R}$, where $\X$ is a two-dimensional real Banach space whose unit sphere is given by:
\begin{align}\label{Hexagon}
S_\X:=\left\{(x,y)\in \mathbb{R}^2:\frac{\sqrt{3}|y|+|x|+\left|\frac{|y|}{\sqrt{3}}-|x|\right|}{2}=1\right\}. 
\end{align}

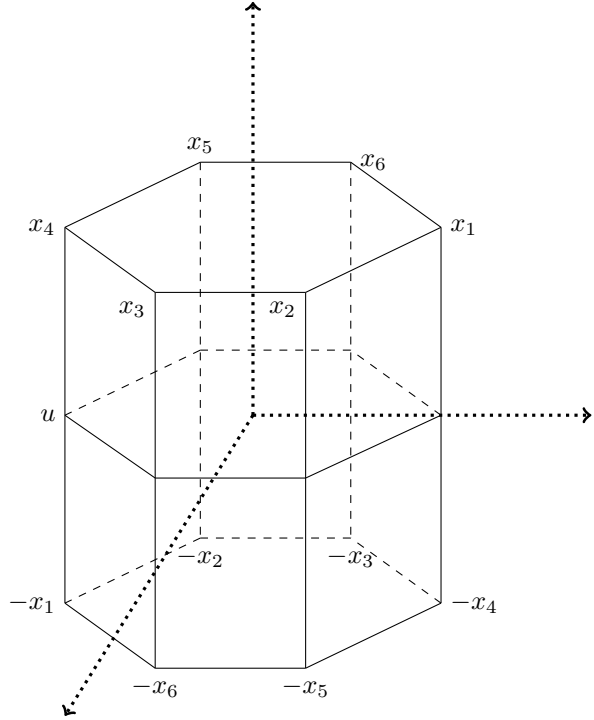
\begin{figure}[h]
\begin{center}
\begin{tikzpicture}

\draw[] (2,0)--(0.8,0.866);
\draw[] (2,0) node[right]{$x_1$};
\draw[] (0.8,0.866)--(-1.2,0.866);
\draw[] (0.8,0.866) node[right]{$x_6$};
\draw[] (-1.2,0.866)--(-3,0);
\draw[] (-1.2,0.866) node[above]{$x_5$};
\draw[] (-3,0)--(-1.8,-0.866);
\draw[] (-3,0) node[left]{$x_4$};
\draw[] (-1.8,-0.866)--(0.2,-0.866);
\draw[] (-1.8,-0.866) node[below left]{$x_3$};
\draw[] (0.2,-0.866)--(2,0);
\draw[] (0.2,-0.866) node[below left]{$x_2$};

\draw[] (2,0)--(2,-5);
\draw[dashed] (0.8,0.866)--(0.8,-4.134);
\draw[dashed] (-1.2,0.866)--(-1.2,-4.134);
\draw[] (-3,0)--(-3,-5);
\draw[] (-1.8,-0.866)--(-1.8,-5.866);
\draw[] (0.2,-0.866)--(0.2,-5.866);

\draw[dashed] (2,-2.5)--(0.8,-1.634);
\draw[dashed] (0.8,-1.634)--(-1.2,-1.634);
\draw[dashed] (-1.2,-1.634)--(-3,-2.5);
\draw[] (-3,-2.5)--(-1.8,-3.336);
\draw[] (-1.8,-3.336)--(0.2,-3.336);
\draw[] (0.2,-3.336)--(2,-2.5);
\draw[] (-3,-2.5) node[left]{$u$};

\draw[dashed] (2,-5)--(0.8,-4.134);
\draw[dashed] (2,-5) node[right]{$-x_4$};
\draw[dashed] (0.8,-4.134)--(-1.2,-4.134);
\draw[dashed] (0.8,-4.134) node[below]{$-x_3$};
\draw[dashed] (-1.2,-4.134)--(-3,-5);
\draw[dashed] (-1.2,-4.134)node[below]{$-x_2$};
\draw[] (-3,-5)--(-1.8,-5.866);
\draw[] (-3,-5) node[left]{$-x_1$};
\draw[] (-1.8,-5.866)--(0.2,-5.866);
\draw[] (-1.8,-5.866) node[below]{$-x_6$};
\draw[] (0.2,-5.866)--(2,-5);
\draw[] (0.2,-5.866) node[below]{$-x_5$};

\draw[very thick,dotted, ->] (-0.5,-2.5)--(4,-2.5);
\draw[very thick, dotted, ->] (-0.5,-2.5)--(-3,-6.5);
\draw[very thick, dotted, ->] (-0.5,-2.5)--(-0.5,3);

\end{tikzpicture}
\caption{Unit sphere of ${\X\oplus_\infty \mathbb{R}}$}
\end{center}
\end{figure}

It is not difficult to see that $S_\X$ is a regular hexagon in $\mathbb{R}^2$ and $S_\mathbb{Z}$ is a hexagonal prism in $\mathbb{R}^3$. Also, $ext(B_\mathbb{Z})=\{\pm x_1,\pm x_2, \pm x_3, \pm x_4,\pm x_5, \pm x_6\}$, where $x_1=(1,0,1)$, $x_2=(\frac{1}{2}, \frac{\sqrt{3}}{2},1)$, $x_3=(-\frac{1}{2}, \frac{\sqrt{3}}{2},1)$, $x_4=(-1,0,1)$, $x_5=(-\frac{1}{2}, -\frac{\sqrt{3}}{2},1)$, $x_6=(\frac{1}{2}, -\frac{\sqrt{3}}{2},1)$. A pictorial description of $S_\mathbb{Z}$ can be seen from Figure - 1.\\

Let $g: \mathbb{Z}\to \mathbb{R}$ be defined by
$$g(x,y,z)=\frac{x+\sqrt{3}y-z}{3}\qquad \forall~ (x,y,z)\in \mathbb{Z}.$$
A simple computation reveals that $|g(x_i)|<1$ for all $i\in \{1,2,3,4,6\}$ and $|g(x_5)|=1.$ Thus, $\|g\|=1$ and $\M_g=\{\pm x_5\}$. Now, define $T:\mathbb{Z}\to \mathbb{Z}$ by
$$T(x,y,z)=g(x,y,z)u\qquad \forall~ (x,y,z)\in \mathbb{Z},$$
where $u=(-1,0,0)$. Our aim is to show that $T$ is not smooth with respect to the usual operator norm but $T$ is smooth with respect to the numerical radius norm.\\

Given any $(x,y,z)\in S_\mathbb{Z}\setminus \{\pm x_5\}$, 
\begin{align*}
\|T(x,y,z)\| & = \|g(x,y,z)u\|\\
& = |g(x,y,z)|\|u\|\\
& = |g(x,y,z)|\\
& < 1.
\end{align*}
On the other hand, $\|T(x_5)\|=1.$ Therefore, $\|T\|=1$ and $\M_T=\{\pm x_5\}.$ Note that $u$ is a non-smooth point and $T(x_5)=-u$. Therefore, it follows from Theorem \ref{PSG-4.1,4.2} that $T$ is not smooth with respect to the usual operator norm.\\

Next, let 
$$\Lambda:=\{\pm (x_5,h): h\in J(x_5)\}.$$
Observe that for any $(v,f)\in \J\setminus \Lambda$,
$$|f(Tv)|=|f(g(v)u)|\leq \|f\|\|g(v)u\|= |g(v)|\|u\|<1.$$
Let
$$f_1(x,y,z)=z, f_2(x,y,z)=-\frac{2}{\sqrt{3}}y,f_3(x,y,z)=-x-\frac{1}{\sqrt{3}}y\qquad \forall~ (x,y,z)\in \mathbb{Z}.$$
Note that $f_1$, $f_2$ and $f_3$ are support functionals of $B_{\mathbb{Z}}$ at $x_5$ and contained in $ext(B_{\mathbb{Z}^*}).$ 
Consequently,
$$J(x_5)=\left\{\lambda_1f_1+\lambda_2f_2+\lambda_3f_3: \lambda_1, \lambda_2, \lambda_3\geq 0,~\lambda_1+\lambda_2+\lambda_3=1\right\}.$$
Evidently, $Tx_5 \in \ker f_1\cap \ker f_2$ and $f_3(Tx_5)=1.$
Therefore, for any $(x_5, h)\in \Lambda,$
$$|h(Tx_5)|\leq 1,$$
and the equality holds for $h=f_3.$ This shows that $\|T\|_w=1,$ and $\M_{\W(T)}=\{(x_5,f_3),(-x_5,-f_3)\}$. Consequently, $T$ is nu-smooth by Theorem \ref{Characterization of NU-smooth in R}.
\end{example}

\begin{example}\label{Operator smooth but not nu-smooth}
Let $\X$ be a two-dimensional Banach space whose unit sphere is a regular hexagon in $\mathbb{R}^2$ defined in Example \ref{NU smooth but not operator smooth}.
\begin{figure}[h]
\begin{center}
\begin{tikzpicture}
\draw[very thick, dotted, ->] (0,0)--(4,0);
\draw[very thick, dotted, ->] (0,0)--(0,4);

\draw[] (0,1.866) node[above right]{$u$};
\draw[] (3,0)--(1.8,1.866);
\draw[] (3,0) node[below right]{$x_1$};
\draw[] (1.8,1.866)--(-1.8,1.866);
\draw[] (1.8,1.866) node[right]{$x_2$};
\draw[] (-1.8,1.866)--(-3,0);
\draw[] (-1.8,1.866) node[above]{$x_3$};
\draw[] (-3,0)--(-1.8,-1.866);
\draw[] (-3,0) node[left]{$-x_1$};
\draw[] (-1.8,-1.866)--(1.8,-1.866);
\draw[] (-1.8,-1.866) node[below left]{$-x_2$};
\draw[] (1.8,-1.866)--(3,0);
\draw[] (1.8,-1.866) node[below left]{$-x_3$};
\end{tikzpicture}
\caption{Unit sphere of $\X$}
\end{center}
\end{figure}
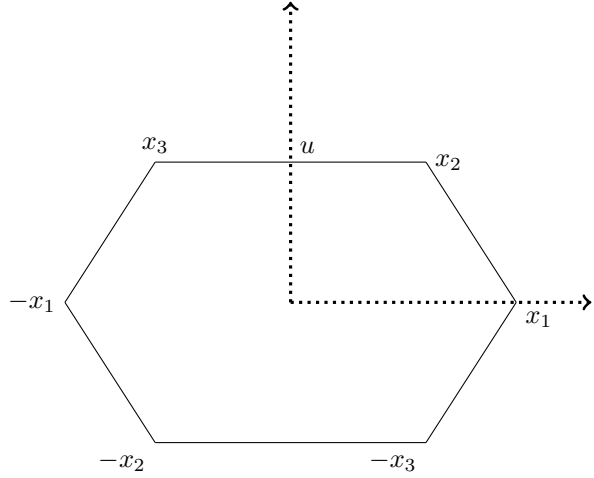
Clearly, $ext(B_\X)=\{\pm x_1,\pm x_2, \pm x_3\},$ where $x_1=(1,0),$ $x_2=(\frac{1}{2}, \frac{\sqrt{3}}{2}),$ $x_3=(-\frac{1}{2}, \frac{\sqrt{3}}{2})$. A pictorial description of $S_\X$ can be seen from Figure - 2.\\ 

Evidently, $ext(B_{\X^*})=\{\pm f_1, \pm f_2, \pm f_3\},$ where 
$$f_1(x,y)=x-\frac{1}{\sqrt{3}}y, f_2(x,y)=x+\frac{1}{\sqrt{3}}y, f_3(x,y)=\frac{2}{\sqrt{3}}y  \qquad \forall~ (x,y)\in \X.$$

Let $g:\X\to \mathbb{R}$ be defined by
$$g(x,y)=x\qquad \forall~ (x,y)\in\X.$$
Define $T:\X \to \X$ by
$$T(x,y)=g(x,y)u\qquad \forall~ (x,y)\in \X,$$ 
where $u=(0,\frac{\sqrt{3}}{2}).$ Our aim is to show that $T$ is not nu-smooth but $T$ is smooth with respect to the usual operator norm.\\

Clearly, $\M_T=\M_g=\{\pm (1,0)\}$ and $T(1,0)=u=(0,\frac{\sqrt{3}}{2})$ is a smooth point of $S_\X$. Therefore, it follows from Theorem \ref{PSG-4.1,4.2} that $T$ is smooth with respect to the usual operator norm. Observe that
$J(x_1)\cap ext(B_{\X^*})=\{f_1,f_2\},$ $J(x_2)\cap ext(B_{\X^*})=\{f_2,f_3\}$ and $J(x_3)\cap ext(B_{\X^*})=\{f_3,-f_1\}$. Now,
\begin{align*}
& f_1(Tx_1)=f_1(g(x_1)u)=f_1(u)=-\frac{1}{2},\\ 
& f_2(Tx_1)=f_2(g(x_1)u)=f_2(u)=\frac{1}{2};\\
& f_2(Tx_2)=f_2(g(x_2)u)=\frac{1}{2}f_2(u)=\frac{1}{4},\\
& f_3(Tx_2)=f_3(g(x_2)u)=\frac{1}{2}f_3(u)=\frac{1}{2};\\
& f_3(Tx_3)=f_3(g(x_3)u)=-\frac{1}{2}f_3(u)=-\frac{1}{2},\\
-& f_1(Tx_3)=-f_1(g(x_3)u)=\frac{1}{2}f_1(u)=\frac{1}{4}.
\end{align*}
Thus, it follows from (\ref{Polygonal reformulation}) that $\|T\|_w=\frac{1}{2}$ and $\{(x_1,f_1),(-x_1,-f_1),(x_1,f_2),(-x_1,-f_2)\}$ is a subset of $ \M_{\W(T)}.$ Consequently, $T$ is not nu-smooth by Theorem \ref{Characterization of NU-smooth in R}.
\end{example}

\end{document}